\documentclass[11pt]{article}

\usepackage{amsfonts,epsfig,fullpage}
\usepackage{caption}
\usepackage{graphicx}
\usepackage{amsmath}
\usepackage[table]{xcolor}

\newtheorem{theo}{Theorem}

\newtheorem{lem}{Lemma}
\newtheorem{claim}[lem]{Claim}
\newtheorem{coro}[lem]{Corollary}

\newtheorem{define}{Definition}

\newtheorem{alg}{Algorithm}
\newtheorem{proc}{Procedure}

\newcommand{\BE}{\begin{enumerate}} \newcommand{\EE}{\end{enumerate}}
\newcommand{\BI}{\begin{itemize}} \newcommand{\EI}{\end{itemize}}
\newcommand{\BDes}{\begin{description}}\newcommand{\EDes}{\end{description}
}
\newcommand{\BT}{\begin{theo}} \newcommand{\ET}{\end{theo}}
\newcommand{\BL}{\begin{lem}} \newcommand{\EL}{\end{lem}}
\newcommand{\BD}{\begin{define}} \newcommand{\ED}{\end{define}}
\newcommand{\BCM}{\begin{claim}} \newcommand{\ECM}{\end{claim}}
\newcommand{\BC}{\begin{coro}} \newcommand{\EC}{\end{coro}}
\newcommand{\BA}{\begin{alg}} \newcommand{\EA}{\end{alg}}
\newcommand{\BP}{\begin{proc}} \newcommand{\EP}{\end{proc}}

\def\FullBox{\hbox{\vrule width 8pt height 8pt depth 0pt}}
\newcommand{\qed}{\;\;\;\FullBox}
\newenvironment{proof}{\noindent{\bf Proof:~~}}{\(\qed\)}
\newcommand{\BPF}{\begin{proof}} \newcommand {\EPF}{\end{proof}}
\newenvironment{proofof}[1]{\noindent{\bf Proof of {#1}:~~}}{\(\qed\)}
\newcommand{\BPFOF}{\begin{proofof}} \newcommand {\EPFOF}{\end{proofof}}

\newcommand{\BEQ}{\begin{equation}} \newcommand{\EEQ}{\end{equation}}
\newcommand{\BEQN}{\begin{eqnarray}}\newcommand{\EEQN}{\end{eqnarray}}
\newcommand{\eqdef}{\stackrel{\rm def}{=}}

\newlength{\saveparindent}
\newlength{\saveparskip}
\providecommand{\keywords}[1]{\textbf{\textit{Key Words---}} #1}

\begin{document}

\title{Circulant almost cross intersecting families }
\author{
    Michal Parnas \\
    The Academic College\\
    of Tel-Aviv-Yaffo \\
     {\tt michalp@mta.ac.il}
}

\maketitle

\abstract{
Let $\mathcal{F}$ and $\mathcal{G}$ be two $t$-uniform families of  subsets over $[k] \eqdef \{1,2,...,k\}$, where $|\mathcal{F}| = |\mathcal{G}|$,
and let $C$ be the adjacency matrix of the bipartite graph whose vertices are the subsets in $\mathcal{F}$ and $\mathcal{G}$,
and there is an edge between $A\in \mathcal{F}$ and $B \in \mathcal{G}$ if and only if $A \cap B \neq \emptyset$.
The pair $(\mathcal{F},\mathcal{G})$ is $q$-almost cross intersecting if every row and column of $C$ has exactly $q$ zeros.

We consider $q$-almost cross intersecting pairs that have a circulant intersection matrix  $C_{p,q}$,
determined by a column vector with $p > 0$ ones followed by $q > 0$ zeros.
This family of matrices includes the identity matrix in one extreme, and the adjacency matrix of the bipartite crown graph in the other extreme.

We give constructions of pairs $(\mathcal{F},\mathcal{G})$ whose intersection matrix is $C_{p,q}$,
for a wide range of values of the parameters $p$ and $q$, and in some cases also prove matching upper bounds.
Specifically, we prove results for the following values of the parameters:
(1) $1 \leq p \leq 2t-1$ and  $1 \leq q \leq k-2t+1$.
(2)  $2t \leq p \leq t^2$ and any $q> 0$, where $k \geq p+q$.
(3)  $p$ that is exponential in $t$, for large enough $k$.

Using the first result we show that if   $k \geq 4t-3$ then
$C_{2t-1,k-2t+1}$ is a maximal isolation submatrix of size $k\times k$ in the $0,1$-matrix $A_{k,t}$,
whose rows and columns are labeled by all  subsets of size $t$ of $[k]$,
and there is a one in the entry on row $x$ and column $y$ if and only if subsets $x,y$ intersect.
}\\



\section{Introduction}
One of the fundamental results of extremal set theory is the theorem of Erd\H os, Ko and Rado~\cite{erdos1961intersection},
which shows that the size of an intersecting $t$-uniform family of subsets over $[k] \eqdef \{1,2,...,k\}$ is bounded above by ${k-1 \choose t-1}$.
Numerous variations of the original problem have been suggested and studied over the years, and among them is the
problem of cross intersecting families of subsets. Specifically, if
$\mathcal{F}$ and $\mathcal{G}$ are two $t$-uniform families of subsets over $[k]$, then
the pair $(\mathcal{F},\mathcal{G})$ is cross intersecting if every subset in $\mathcal{F}$ intersects with every subset in $\mathcal{G}$ and vice versa.
Pyber~\cite{pyber1986new} proved that in this case $|\mathcal{F}|\cdot|\mathcal{G}| \leq {k-1 \choose t-1}^2$.

Many of the extremal combinatorial problems considered so far can be inferred as results about maximal submatrices of the $0,1$-matrix $A_{k,t}$
of size ${k \choose t}\times{k \choose t}$, whose rows and columns are labeled by all  subsets of size $t$ of $[k]$,
and there is a one in the entry on row $x$ and column $y$ if and only if subsets $x,y$ intersect.
Hence, in this setting, the result of Erd\H os, Ko and Rado can be inferred as stating the size of the largest all-one square principal submatrix  of $A_{k,t}$,
and the result of Pyber states the size of the largest all-ones submatrix of $A_{k,t}$.
We note that considering the classic results of extremal combinatorics as maximal submatrices of $A_{k,t}$,
allows us to employ tools from algebra in addition to the combinatorial techniques.

Another variation of the problem of cross intersecting families was introduced by Gerbner et al.~\cite{gerbner2013almost},
which defined the notion of a $q$-almost cross intersecting pair $(\mathcal{F},\mathcal{G})$.
Here every subset in $\mathcal{F}$  does not intersect with exactly $q$ subsets in $\mathcal{G}$ and vice versa.
If  $\mathcal{F} = \mathcal{G}$ then $\mathcal{F}$ is called a $q$-almost intersecting family of subsets.
Hence, if $C$ is the submatrix of $A_{k,t}$ whose rows are indexed by the subsets of $\mathcal{F}$ and columns by the subsets of $\mathcal{G}$,
then every row and column of $C$ has exactly $q$ zeros.
Using a classic theorem of Bollob{\'a}s~\cite{bollobas1965generalized}, it is possible
to prove that the largest square submatrix $C$ of $A_{k,t}$, representing a $1$-almost cross intersecting pair,  is of size ${2t \choose t}\times {2t \choose t}$.
A theorem proved by~\cite{gerbner2013almost}, shows that if $C$ is a submatrix of size $n \times n$ of $A_{k,t}$, with exactly $q$ zeros in
each row and column, then $n \leq (2q - 1){2t \choose t}$.

In this paper we consider the problem of finding maximal {\em circulant} submatrices of $A_{k,t}$,
representing an almost cross intersecting pair, for a range of parameters.
A circulant matrix is a matrix in which each row is shifted one position to the right compared to the preceding row
(or alternatively, each column is shifted one position compared to the preceding column).
Therefore, such a matrix $C$ is determined completely by its first row or first column.
Circulant matrices were studied extensively in the context of the multiplicative commutative semi-group of circulant Boolean matrices
(see for example~\cite{butler1974circulant, schwarz1974circulant, daode1990circulant, butler1976semigroup}).
However, they were not studied in the context of extremal combinatorics, besides some special cases that will be discussed shortly.

\begin{figure}[htb!]
\captionsetup{width=0.8\textwidth}
	$$
C_{4,4} = \left(
\begin{array}{cccccccc}
 1 & 0 & 0 & 0 & 0  & 1 & 1 & 1    \\
 1 & 1 & 0 & 0  & 0 & 0 & 1 & 1 \\
 1 & 1 & 1 & 0  & 0 & 0 & 0 & 1    \\
 1 & 1 & 1 & 1  & 0 & 0 & 0 & 0    \\
 0 & 1 & 1 & 1 & 1  & 0 & 0 & 0 \\
 0 & 0 & 1 & 1 & 1  & 1 & 0 & 0   \\
 0 & 0 & 0 & 1 & 1  & 1 & 1 & 0     \\
 0 & 0 & 0 & 0 & 1  & 1 & 1 & 1    \\
\end{array}
\right)
	$$
\caption{The circulant matrix $C_{p,q}$, where $p = q = 4$.}
\label{fig:C44}	
\end{figure}

Our focus will be on circulant matrices that are determined by a column vector with $p$ ones followed by $q$ zeros.
Such a matrix will be denoted by $C_{p,q}$. See  Figure~\ref{fig:C44}	for an example.
Thus, in one extreme, if $p=1$ and $q > 0$, then $C_{p,q}$ is the identity matrix.
The other extreme  is  $q = 1$ and $p > 0$, and then  $C_{p,q}$ is the adjacency matrix of a crown graph
(where a crown graph is a complete bipartite graph from which the edges of a perfect matching have been removed).
Hence, the structure of the circulant matrix $C_{p,q}$  forms a bridge which connects these two extreme cases,
and it is interesting to find a unifying theorem
which determines the maximal size of the matrix $C_{p,q}$ as a function of $p,q,k$ and $t$.

We note that two trivial examples of circulant submatrices of $A_{k,t}$ include the case of $q= 0$, where we get
an all-one submatrix of $A_{k,t}$ of maximal size  ${k-1 \choose t-1} \times {k-1 \choose t-1}$,
and the case of $p = 0$, where we get an  all-zero submatrix of $A_{k,t}$ of maximal size ${k/2 \choose t}\times {k/2 \choose t}$.
Hence, the problem of studying the  size of circulant submatrices $C_{p,q}$ of $A_{k,t}$
is interesting only if both $p,q > 0$.
Furthermore, we must require that $k \geq 2t$, as otherwise, $A_{k,t}$ is the all-one matrix itself.

As we discuss shortly, one of our results also provides a simple construction of maximal isolation submatrices of $A_{k,t}$,
thus providing simple small witnesses to the Boolean rank of $A_{k,t}$.
The  {\em Boolean rank} of a $0,1$-matrix $A$ of size $n \times m$ is equal to the smallest integer $r$,
such that $A$ can be factorized as a product of two $0,1$-matrices, $X\cdot Y=A$,   where
$X$ is a matrix of size $n \times r$  and $Y$ is a matrix of size $r \times m$, and all additions and multiplications are Boolean
(that is, $1+1= 1$, $1+0 = 0+1 = 1$, $1 \cdot 1 = 1, 1 \cdot 0 = 0 \cdot 1 = 0$).
A $0,1$-matrix $B$ of size $s\times s$ is called an {\em isolation matrix}, if it contains $s$ ones,
such that no two of these ones  are in the same row or column of $B$,
and no two of these ones are contained in a $2\times 2$ all-one submatrix of $B$.
It is well known that if $B$ is an isolation submatrix of size $s\times s$ in a given $0,1$-matrix $A$,
then $s$ bounds below the Boolean rank of $A$ (see~\cite{ beasley2012isolation, KN97}).
However, computing the Boolean rank or finding a maximal isolation submatrix in general is an NP-hard problem
(see~\cite{  gruber2007inapproximability,   orlin1977contentment,shitov2013complexity}).
Hence, it is interesting to find and characterize families of maximal isolation sets for specific given matrices.


\subsection{Our Results}
Our main goal is to determine the range of parameters, $p$ and $q$, for which  $C_{p,q}$ is a submatrix of $A_{k,t}$.
The constructions and upper bounds we present differ in their structure and proof methods according to the
size of $p,q$ compared to $t,k$.

We first consider the range of values of relatively small $p$, that is $1 \leq p \leq 2t-1$,
and prove in Section~\ref{smallp} the following positive result.\\

{\bf Theorem 1:}
{\it
Let $k \geq 2t$, $1 \leq p \leq 2t-1$ and  $1 \leq q \leq k-2t+1$. Then $C_{p,q}$ is a submatrix of $A_{k,t}$.
}\\

In the extreme case of $p = 1$ and $q = k-2t+1$, this construction gives the identity submatrix of size $(k-2t+2)\times (k-2t+2)$.
Recently,~\cite{parnas2019maximal} proved that this is the maximal size of an identity submatrix in $A_{k,t}$.

The other extreme is  $p = 2t-1$ and $q = k-2t+1$, in which case we get a circulant submatrix of size $k\times k$.
As we show in Section~\ref{smallp}, if   $k \geq 4t-3$ then
$C_{2t-1,k-2t+1}$ is a maximal  isolation submatrix of size $k\times k$ in $A_{k,t}$.
Since the Boolean rank of $A_{k,t}$ is $k$ for $k \geq 2t$ (see~\cite{parnas2019boolean}),
then the size of a maximal isolation submatrix of $A_{k,t}$ is upper bounded by $k\times k$,
and thus, our result is optimal in this case.

Furthermore, for $k = 2t + p -2$ and $p \geq 2$, the construction described in Theorem~1 provides an isolation submatrix of size $(2p-1)\times (2p-1)$.
We note that~\cite{parnas2019maximal} gave constructions of isolation submatrices in $A_{k,t}$, of the same size as achieved here.
However, the constructions described in~\cite{parnas2019maximal} are quite complex, and thus, the result described in Theorem~1
provides an alternative simpler construction of a maximal isolation submatrix in $A_{k,t}$, for large enough $k$.

We then prove the following upper bound that matches the size of the construction given in Theorem 1, for the range of values of $1 \leq p \leq 2t-1$
and $q \geq p-1$. The proof of this result characterizes the structure of the Boolean decompositions of $C_{p,q}$ for this range of parameters. \\

{\bf Theorem 2:}
{\it
Let $C_{p,q}$ be a submatrix of $A_{k,t}$, where  $k \geq 2t$, $ 1 \leq p \leq 2t-1$ and $q > 0$.
If $q \geq p-1$ then   $q \leq k-2t + 1$.
}\\

In Section~\ref{midp} we address the range of slightly larger values of $p$, that is, $2t \leq p \leq t^2$,
and provide a different construction of circulant submatrices of $A_{k,t}$ of the form $C_{p,q}$.
As we show, for this range of values of $p$, there is no upper bound on the size of $q$, as we had in Theorems~1 and~2,
as long as $k  \geq p+q$.

Furthermore, the proof for this range of parameters provides a decomposition of $C_{p,q}$ into a product of two Boolean circulant matrices $X,Y$,
where $X$ has $t$ ones in each row and $Y$ has $t$ ones in each column.
If we view the rows of $X$ and the columns of $Y$ as the characteristic vectors of subsets of size $t$,
then $X$ and $Y$ each represents a {\em circulant $t$-uniform family}.
Thus, the construction used in the proof of the next theorem,
uses a pair $\mathcal{F}, \mathcal{G}$ of circulant families to construct $C_{p,q}$.\\
%

{\bf Theorem 3:}
{\it
Let $2t \leq p \leq t^2$ and $q > 0$.  Then $C_{p,q}$ is a submatrix of  $A_{k,t}$ for $k \geq q+p$.
}\\

Finally, in Section~\ref{seclargep}, we consider the range of large $p$.
Using the result of~\cite{gerbner2013almost} stated above,  we know that if $C_{p,q}$ is a submatrix of $A_{k,t}$ of size $n\times n$,
then $n \leq (2q - 1){2t \choose t}$, and~\cite{scott2014hypergraphs} proved a conjecture
of~\cite{gerbner2012almost} and showed that for a large enough $q$ and $t$, the size of a
$q$-almost intersecting family $\mathcal{F}$  is bounded by $(q+1){2t-2 \choose t-1}$.
Note that this last result refers to $q$-almost cross intersecting pairs $(\mathcal{F},\mathcal{G})$ in which $\mathcal{F} =\mathcal{G}$.
Furthermore, the constructions presented in~\cite{scott2014hypergraphs}, which achieve this bound, do not have a circulant intersection matrix.

Indeed,  we can get a better upper bound for circulant submatrices of the form $C_{p,q}$.
Using a theorem of Frankl~\cite{frankl1982extremal} and Kalai~\cite{kalai1984intersection} about skew matrices,
it is possible to show that $p \leq {2t \choose t}-1$. Hence,  if $C_{p,q}$ is a submatrix of size $n\times n$ of $A_{k,t}$ then  $n\leq {2t \choose t} + q-1$.

In the extreme case of $p = {2t \choose t}-1$ and $q =1$,
the simple construction that takes all  subsets of size $t$ of $[2t]$ as row and column indices,
results in a submatrix $C_{p,q}$ of size ${2t \choose t}\times {2t \choose t}$.
This is optimal, as it matches the upper bound of ${2t \choose t} + q-1$.

For larger $q$, we give a simple construction of $C_{p,q}$ for $p = q\cdot ({2t/q \choose t/q} -1) $,
when $t \bmod q = 0$ and $k$ is large enough.
Note that there is a relatively large gap between the size of $C_{p,q}$ stated here, and the upper bound of ${2t \choose t} + q-1$.
As we prove, this gap can be slightly narrowed for $q = 2$:\\


{\bf Theorem 4:}
{\it
Let $q = 2$ and $p = 2^t + 2^{t-2} -2$, where $t > 2$. Then $C_{p,q}$ is a submatrix of $A_{k,t}$ for large enough $k$.
}\\

We conclude by  considering the case of $t = 2$ and $p = {2t \choose t}-1 = 5$ and fully characterize it.
As we show, in this case, $C_{p,q}$ is a submatrix of $A_{k,t}$, for $q= 1$ and $k \geq 5$, or for $q = 3$ and $k \geq 6$.
Thus, for $t = 2, p = 5$ and $q = 1,3$, we get a result which  matches the upper bound of ${2t \choose t} + q-1$.
However, as we prove, for $t = 2, p = 5$ and $q > 0$, $q \neq 1,3$, there is no $k$ for which $C_{p,q}$ is a submatrix of $A_{k,t}$.
This implies that the upper bound of ${2t \choose t} + q-1$ is not tight in general.
It remains an open problem to determine for what values of $q > 1$ is $C_{p,q}$ a submatrix of $A_{k,t}$,
given that $ p = {2t \choose t}-1$ and $t > 2.$




\section{The range of $1 \leq p \leq 2t-1$}
\label{smallp}


In this section we prove  Theorems 1 and 2, which address the range of small $p$, that is, $1 \leq p \leq 2t-1$.
As stated above, this range of values includes the identity matrix, as well as allows us to provide a simple construction of maximal isolation sets
for large enough $k$.

It will be useful to  identify subsets of $[k]$  with their characteristic vectors.
Thus, a  subset of size $t$ of $[k]$ will be represented by a $0,1$-vector of length $k$ with exactly $t$ ones.
Furthermore, in order to show that some matrix $C$ of size $n \times m$ is a submatrix of $A_{k,t}$,
it will be enough to show that there exists a Boolean decomposition  $C = X \cdot Y $, where
$X$ is a Boolean matrix of size $n \times k$ with exactly $t$ ones in each row,
and $Y$ is a Boolean matrix of size $k \times m$ with exactly $t$ ones in each column,
and all operations are Boolean.

\subsection{A construction of $C_{p,q}$ for $1 \leq p \leq 2t-1$}

The following lemma will be useful in proving Theorem~1. It shows that it is possible to  decompose a matrix of the form $C_{p,q}$
into a product of two circulant matrices of the same type, for a wide range of parameters.

\BL
\label{productcirculant}
Let $i,j,z$ be three integers, such that $i , j \geq 1$ and $i+j-1 \leq z$.
Then $$ C_{i,z-i} \cdot C_{j,z-j} = C_{i+j-1,z-i-j+1} .$$
\EL
\BPF
It is well known that the product of two circulant Boolean matrices is a circulant Boolaen matrix (where all operations are Boolean).
Thus, it is enough to determine the first column $c = (c_1,c_2,...,c_z)$
of the product matrix $C_{i,z-i} \cdot C_{j,z-j} $, and to show that it has $i+j-1$ ones, followed by $z-i-j+1$  zeros.

First it is clear that $c_s = 1$ for $1 \leq s \leq i$,
since the first element in each of the first $i$ rows of $C_{i,z-i}$ is a $1$,
and the first element of the first column of $C_{j,z-j} $ is also a $1$ (since $i,j \geq 1$).

Next consider element $c_{i+s}$ for  $1 \leq s \leq j-1$.
Note that row $i+s$  of $C_{i,z-i}$ begins with $s$ zeros and then has $i$ ones,
and  the first $j$ elements of the first column of $ C_{j,z-j} $ are ones.
Since $s \leq j-1$, then the result of multiplying row $i+s$ of $C_{i,z-i}$ with the first column of $ C_{j,z-j} $, is a one.

Hence, we showed so far that the first $i+j-1$ elements of the vector $c$ are ones, and it remains to show that the remaining elements of $c$ are all zeros.
But row $i+j$ of $C_{i,z-i}$ begins with $j $ zeros. Therefore, multiplying it with the first column of $ C_{j,z-j} $, results in a zero.
The remaining rows of $C_{i,z-i}$ begin with at least $j $ zeros, and so also contribute a zero when multiplied with the first column of $ C_{j,z-j} $.
\EPF\\

Using Lemma~\ref{productcirculant}, we can now prove Theorem 1.

\BT
\label{circulargeneral}
Let $ k \geq 2t$, $1 \leq p \leq 2t-1$ and $1 \leq q \leq k-2t+1$. Then $C_{p,q}$ is a  submatrix of $A_{k,t}$.
\ET
\BPF
Let $1 \leq i,j \leq t$ such that $i + j-1 = p$. Thus, $1 \leq p \leq 2t-1$, as claimed.
Let $J_{n,m}$ be the all-one matrix of size $n \times m$, and $O_{n,m}$ the all-zero matrix of size $n \times m$.
Define, two matrices $X$ and $Y$ as follows:
$$X = [C_{i,p+q-i} O_{p+q,t-j} J_{p + q,t-i}], \; \; \; \; \; \;
Y = \left[
  \begin{array}{c}
    C_{j,p+q-j} \\
    J_{t-j,p+q}\\
    0_{t-i,p+q} \\
  \end{array}
\right].
$$
Using Lemma~\ref{productcirculant}, where $z = p+q$, we have that
$$X\cdot Y = C_{i,p+q-i} \cdot C_{j,p+q-j} =  C_{i+j-1,p+q-i-j+1} =  C_{p,q}.$$
Furthermore, each row of $X$ and each column of $Y$ is a vector with exactly $t$ ones, whose length is:
$$(p+q) + (t-j) + (t-i) = p+q + 2t  - i - j =  p+q + 2t - (p+1) = q + 2t - 1.$$
Therefore, if $k \geq q+2t-1$, then we can view the rows of $X$ and columns of $Y$ as
the characteristic vectors of subsets in ${[k] \choose t}$.
Thus, $X\cdot Y = C_{p,q}$ is a  submatrix of $A_{k,t}$ as claimed.
\EPF\\

As we show next, if   $k \geq 4t-3$ then the construction described in the proof of Theorem~\ref{circulargeneral},
provides a maximal isolation submatrix of size $k\times k$ in $A_{k,t}$. This result is optimal since
the Boolean rank of $A_{k,t}$ is $k$ for $k \geq 2t$ (see~\cite{parnas2019boolean}).

\begin{coro}
Let $2 \leq p \leq 2t-1$ and let $k = 2t+p-2$. Then $C_{p,p-1}$ is an isolation submatrix of size $(2p-1)\times (2p-1)$ in $A_{k,t}$.
Furthermore, if $k \geq 4t-3$ then  $C_{2t-1,k-2t+1}$ is an isolation submatrix of  size $k\times k$ in $A_{k,t}$.
\end{coro}
\BPF
Let $k = 2t+p-2$. If we set $q = k-2t + 1 = (2t + p-2) - 2t+ 1 = p-1$,
then by Theorem~\ref{circulargeneral}, $C_{p,q}$ is a submatrix of $A_{k,t}$ of size $(2p-1) \times (2p-1)$ since $p + q = 2p-1$.
It is easy to verify that in this case, since $q = p-1$, then the ones on the main diagonal of $C_{p,q}$ form an isolation set of size $p+q$.

In the extreme case of $ p = 2t-1$, and if $k \geq 4t-3$, then $q = k-2t+1 \geq 2t - 2 \geq p-1$, and we get
an isolation matrix $C_{p,q}$ of size $k \times k$, since $p + k-2t + 1 = k$.
\EPF\\

\subsection{Upper Bounds on the size of $C_{p,q}$ for $1 \leq p \leq 2t-1$}
We now turn to prove Theorem~2, which provides a matching upper bound to the size of the construction given in Theorem~1,
for  $1 \leq p \leq 2t-1$ and $q \geq p-1$.
We note that if $q \geq p-1$ then $p+q \leq k$ (for any value of $p$), since in this case
$C_{p,q}$ is an isolation submatrix of $A_{k,t}$. Thus, its Boolean rank, which is $p+q$, is bounded above by $k$, which is
the Boolean rank of $A_{k,t}$.
However, the proof of Theorem~2, which provides a tight upper bound on $p+q$, will require a more elaborate proof.

The following simple claim is easy to verify, and will be needed for the proof of Theorem~2.

\BCM
\label{clm:monosub}
Let $B$ be an all-one submatrix of size $i\times j$ of $C_{p,q}$, where $p,q > 0$.
Then, $1 \leq i,j \leq p$ and $i+j \leq p+1$.
\ECM

The next lemma is a generalization of a claim proved in~\cite{parnas2019boolean}, which characterizes the Boolean decompositions of the identity matrix.
Here we characterize the Boolean decompositions of circulant isolation matrices of the form $C_{p,q}$.

Denote by $|x|$ the number of ones in a vector $x$,
and let $x \otimes y$ denote the outer product of a column vector $x$ and a row vector $y$, where both $x,y$ are of length $n$.
That is, $x \otimes y$ is a matrix of size  $n\times n$.

\BL
\label{clm:circulantdecompos}
Let $p,q > 0$ and  $n = p+q$.
Let $X\cdot Y = C_{p,q}$ be a Boolean decomposition of  $C_{p,q}$, where
 $X$ is an $n \times r$ Boolean matrix and $Y$ is an $r \times n$ Boolean matrix.
Denote by
$x_1, \ldots, x_r$ the columns of $X$, and by $y_1,\ldots,y_r$ the rows of $Y$.
Then:
\begin{enumerate}
\item
For each $i \in [r]$, if $x_i$ has more than $p$ ones then $y_i$ is the all-zero vector,
and if $y_i$ has more than $p$ ones then $x_i$ is the all-zero vector.
\item
For each $i \in [r]$, if $|x_i|, |y_i| > 0$, then $|x_i| + |y_i| \leq p+1$.
\item
If $q \geq p-1$, then there exist $n $ indices $i_1,...,i_n$, such that $|x_{i_j}|,|y_{i_j}| > 0$ for every $j \in [n]$.
\end{enumerate}
\EL
\BPF
For simplicity, denote $C = C_{p,q}$.
If we write the decomposition $X\cdot Y = C$ with outer products, then $C = \sum_{i=1}^r x_i \otimes y_i$.

First note that if we have an $i$ such that $x_i$ has more than $p$ ones, and $y_i$ is not the all-zero vector, then
$x_i \otimes y_i$ has a column with more than $p$ ones.
Since the addition is the Boolean addition, then $\sum_{i=1}^r x_i \otimes y_i \neq C$.
A similar argument shows that if $y_i$ has more than $p$ ones then $x_i$ is the all-zero vector. Thus, item (1) follows.

Assume now, by contradiction, that item (2) does not hold.
Thus, there exists an $i$, such that $|x_i|, |y_i| > 0$ and $|x_i| + |y_i| > p+1$.
Let $|x_i| = s$ and $|y_i| = \ell$, where by our assumption $\ell \geq p - s + 2$.
Thus, the matrix $x_i \otimes y_i$ has an all-one submatrix $B$ of size $ s  \times \ell$.
Since the addition is Boolean, $C_{p,q}$, also has an all-one submatrix of size $ s \times \ell \geq s \times (p-s+2)$,
in contradiction to Claim~\ref{clm:monosub}.

It remains to prove item (3).
Since $q \geq p-1$, then $C$ is an isolation matrix. Therefore, its Boolean rank is $n = p+q$.
Assume by contradiction that there are strictly less than $n$ pairs $x_i,y_i$ such that $|x_i|, |y_i| >0$.
Note that if $x_i$ or $y_i$ is the all-zero vector then $x_i \otimes y_i$ is the all-zero matrix.
Thus, we can remove from $X$ any column $x_i$ which is the all-zero vector, and remove the corresponding row $y_i$ from $Y$,
and similarly, remove from $Y$ any row $y_j$ which is the all-zero vector,
and remove the corresponding column $x_j$ from $X$.
We get two new matrices $X',Y'$, such that $X'\cdot Y' = C$, where the size of $X'$ is $n \times \ell$, the size of $Y'$ is $\ell \times n$,
and by our assumption $\ell < n$.
Therefore, the Boolean rank of $C$ is strictly less than $n$, and we get a contradiction.
\EPF

\BL
\label{claim:maxones}
Let $p,q > 0$ and  $q \geq p-1$, and let $n = p+q$.
Let $X\cdot Y = C_{p,q}$ be a Boolean decomposition of  $C_{p,q}$, where $X$ is an $n \times r$ Boolean matrix
and $Y$ is an $r \times n$ Boolean matrix. Then the total number of $1$'s in both $X$ and $Y$ is at most $(p+1)n + (r-n)n$.
\EL
\BPF
Let $x_1, \ldots, x_r$ be the columns of $X$, and $y_1,\ldots,y_r$ the rows of $Y$.
By Lemma~\ref{clm:circulantdecompos}, there exist $n$ indices $i_1,...,i_n$,
such that $|x_{i_j}|,|y_{i_j}| > 0$ for every $j \in [n]$.
Furthermore, for these indices it holds that $|x_{i_j}| + |y_{i_j}| \leq p+1$.
Assume, without loss of generality, that these are indices $1,...,n$.

As to the remaining pairs, $x_i,y_i$, for $n < i \leq r$:
by Lemma~\ref{clm:circulantdecompos}, if $|x_i|,|y_i| >0$ then $|x_i|+|y_i| \leq p+1$,
and if $|x_i| \geq p+1$ then $y_i$ is the all-zero vector,
and similarly if $|y_i| \geq p+1$ then $x_i$ is the all-zero vector.
Thus, $|x_i|+|y_i|$ is maximized when  $x_i$ or $y_i$ is the all-zero vector and the other is the all-one vector,
since in this case $|x_i|+|y_i| = n = p+q \geq p+1$.

Hence, the number of ones in both $X$ and $Y$ is at most $(p+1)n + (r-n)n$.
\EPF\\


\BT
Let $C_{p,q}$ be a submatrix of $A_{k,t}$, where  $k \geq 2t$, $ 1 \leq p \leq 2t-1$ and $q > 0$.
If $q \geq p-1$ then $q \leq k-2t + 1$.
\ET
\BPF
Consider the Boolean decomposition $X\cdot Y = A_{k,t}$,
where $X$ is a matrix of size ${k \choose t} \times k$   and $Y$ is a matrix of size $k \times  {k \choose t}$,
and $X$ and $Y$ each contain  all characteristic vectors of subsets in ${[k] \choose t}$.

Since $C_{p,q}$ is a submatrix of $A_{k,t}$ then there exist two matrices
$X' \subseteq X, Y' \subseteq Y$, such that $X'\cdot Y' = C_{p,q}$.
Notice that $X'$ is an $n \times k$ matrix and $Y'$ is an $k \times n$ matrix, where $n = p+q$,
and  the total number of $1$'s in both $X'$ and $Y'$ is exactly $2n t$.
But, by Lemma~\ref{claim:maxones}, the total number of $1$'s in both $X'$ and $Y'$ is at most $(p+1)n+(k-n)n$.
Thus, $ 2n t \le (p+1)n+(k-n)n$. Hence, $p + q = n \le k-2t+ p +1$, as claimed.
\EPF\\

\section{The range of $2t \leq p \leq t^2$}
\label{midp}

The circulant decomposition given in Lemma~\ref{productcirculant} is not suitable for $p \geq 2t$, since if we take
the decomposition $ C_{i+j-1,z-i-j+1} =C_{i,z-i} \cdot C_{j,z-j} $, and let
$p = i+j-1$ and $p \geq 2t$, then $i+j \geq 2t+1$.
Thus, either $i$ or $j$ are strictly larger than $t$, and therefore, the rows of
$C_{i,z-i}$ or the columns of $C_{j,z-j}$ cannot represent subsets of size $t$ of $[k]$.

\begin{figure}[htb!]
\captionsetup{width=0.9\textwidth}
\centering
    \includegraphics[width=0.9\textwidth]{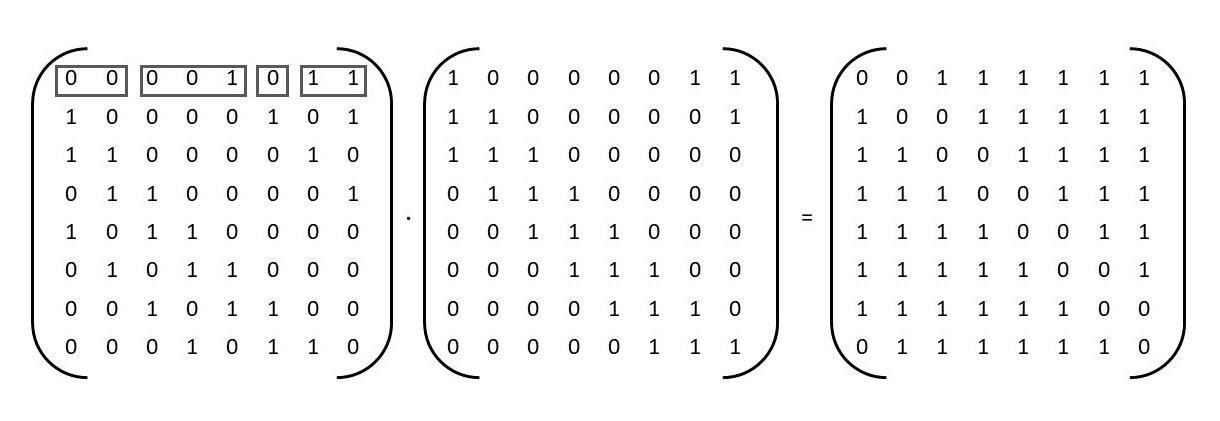}
\caption{\small The construction described in Theorem~3 for $t = 3, p =6, q = 2$,
where the matrices presented are $X \cdot Y = C_{p,q}$.
Here $z = (0,0,1)$, $\ell_1 = \lfloor p/t\rfloor = 2$ and $\ell_2 = p \bmod t = 0$.
The different components of the first row $x_1$ of the matrix $X$ are outlined with rectangles. Row $x_1$ begins with $q$ zeros, followed by the vector $z$,  then
$\ell_1 - 1$ zeros and finally $t-\ell_1 + 1$ ones.}
\label{fig:product}
\end{figure}

However, as we prove next, when $2t \leq p \leq t^2$, there exists a different circulant
decomposition $C_{p,q} = X\cdot Y$,  in which each row of $X$ and each column of $Y$ has exactly $t$ ones as required.
See Figure~\ref{fig:product} for an illustration, and note also that since $2t \leq p \leq t^2$ then $t \geq 2$.

\BT
\label{circularbigp}
Let $2t \leq p \leq t^2$ and $q > 0$. Then $C_{p,q}$ is a  submatrix of  $A_{k,t}$ for $k \geq q+p$.
\ET
\BPF
We prove that there exist two matrices $X,Y$, such that $X \cdot Y = C_{p,q}$,
where $X$ and $Y$ are matrices of size $p+q$, and each one has exactly $t$ ones in each row and column.
It will be  convenient throughout this proof to assume that the first row of $C_{p,q}$ is $(0,0,...,0,1,1,...,1)$, that is, $q$ zeros followed by $p$ ones
(this can be achieved by moving the last row of $C_{p,q}$ as defined in the introduction to be first).

The matrix $Y$ is the circulant matrix defined by the first column  vector with  $t$ consecutive ones followed by $p+q-t$ zeros.
The matrix $X$ is also  circulant and is defined as follows.
Note that each $t$ consecutive columns of $Y$ have a one in a common row. The idea is to spread the $t$ ones in each row of $X$,
such that there are $t-1$ zeros between every two ones of $X$, and each such one of $X$ intersects with $t$ consecutive columns of $Y$.
Since $t$ is not always a divisor of $p$, then the definition of $X$ is slightly more complex.

Let $z = (0,0,...,0,1)$ be a vector of length $t$ that starts with  $t-1$ zeros followed by a one, and let $\ell_1 = \lfloor p/t\rfloor$ and $\ell_2 = p\bmod t$.
Recall also that $(p \bmod t) = p- t\lfloor \frac{p}{t}\rfloor$.
Then $X$ is generated by the following first row vector $x_1$:
\begin{enumerate}
\item
If $\ell_2 = 0$, then:
$$x_1 = (\underbrace{0,0,...,0}_{q},\underbrace{z,z,...,z}_{\ell_1-1},\underbrace{0,0,...,0}_{\ell_1-1}, \underbrace{1,1,...,1}_{t-\ell_1+1 }).$$
Note that $2 \leq \ell_1 \leq t$, since $2t \leq p \leq t^2$. Therefore, $x_1$ is well defined, since
$t- \ell_1  + 1>0$ and $\ell_1 -1 > 0$.
 The length of $x_1$ is $p+q$, as required, since:
$$q + t (\ell_1 -1) + (\ell_1 -1) + (t - \ell_1 +1 )  = q + t \lfloor\frac{p}{t}\rfloor  = q+ p.$$
\item
If $\ell_2 >0$, then:
$$x_1 = (\underbrace{0,0,...,0}_{q},\underbrace{z,z,...,z}_{\ell_1-1}, \underbrace{0,0,...,0}_{\ell_1}, \underbrace{1,1,...,1}_{t-\ell_1},\underbrace{0,0,...,0,1}_{ \ell_2}).$$
In this case $\ell_1 < t$ (since, otherwise, $\ell_2 = 0$), and the length of $x_1$ is again $q+p$:
$$q + t (\ell_1 -1) + \ell_1 + (t - \ell_1 ) + \ell_2  = q + t \lfloor\frac{p}{t}\rfloor + p\bmod t = q+ p.$$
\end{enumerate}
Also note that in both cases, $x_1$ has exactly $t$ ones as claimed.
Finally, we must show that $X\cdot Y = C_{p,q}$.
Since both $X$ and $Y$ are circulant, then the resulting product matrix $W = X\cdot Y$ is circulant.
Thus, it is sufficient to prove that the first row $w = (w_1,w_2,...,w_{p+q})$ of $W$ is equal to the first row
$c = (0,0,...,0,1,1,...,1)$ of $C_{p,q}$. Let $y_1,...,y_{p+q}$ be the columns of $Y$.

We first show that the first $q+(\ell_1 -1)t$ coordinates of $w$ and of $c$ are identical.
Note that for this range of indices, the coordinates of $x_1$ are identical in both definitions of $x_1$ above.
\begin{itemize}
\item
 $w_1,...,w_{q} = 0$:
Columns $y_1,...,y_{q}$ have ones only in positions  at most $q+t-1$, whereas $x_1$ has zeros in the first $q+t-1$ positions.
Hence, $w_i = x_1 \cdot y_i = 0$, for $1 \leq i \leq q$.
\item
$w_{q+i}  = 1$, for $1 \leq i \leq (\ell_1 -1)t$:
By the definition of $x_1$ and $z$, the vector $x_1$ has ones in positions $q+jt$, for $1 \leq j \leq \ell_1 -1$
(this holds for both cases of the definition of $x_1$ above). Furthermore,
column $y_{q+i}$ has ones in positions $q+i,...,q+i+t-1$, since
$$y_{q+i} = (\underbrace{0,0,...,0}_{q+i-1},\underbrace{1,1,...,1}_t,0,0,...,0).$$
Hence, $x_1$ intersects with column $y_{q+i}$, for $(j-1)t+1 \leq i \leq jt$, where $1 \leq j \leq \ell_1 -1$.
This implies that $w_{q+i} = x_1 \cdot y_{q+i} = 1$, for $1 \leq i \leq (\ell_1 -1)t$.
\end{itemize}

We now show that the remaining coordinates of $w$ and $c$ are identical. That is, for coordinates
$q+(\ell_1 -1)t+1 \leq j \leq p+q$, it holds that $w_j = c_j = 1$.
Here there is a different analysis for the two definitions of $x_1$,
which depends on the value of $\ell_2 = p \bmod t$.
\begin{itemize}
\item
$\ell_2 = 0$:
In this case, the vector $x_1$ has a block of $t-\ell_1+1$ consecutive ones starting from position
$q+t(\ell_1-1)+\ell_1 \leq q + t\ell_1$ and until position $p+q$, where the inequality follows since $p \leq t^2$ and therefore, $\ell_1 \leq t$.

The first $y_{q+i}$ we have to consider is $y_{q+i}$ for $i=(\ell_1 -1)t+1$. As before, it has ones in positions $q+i,...,q+i+t-1$.
Thus, it has a one in position $q+i+t-1 = q+ t(\ell_1 -1)+1 + t-1  = q + t\ell_1$, and so it intersects with this block of consecutive ones in $x_1$.
It is also clear that the remaining columns of $Y$ all intersect with this consecutive block of ones of $x_1$,
since they all have ones in positions greater than $q + t\ell_1$.
\item
$\ell_2 > 0$:
In this case,  $x_1$ has a block of $t-\ell_1$ consecutive ones starting from position
$q+t(\ell_1-1)+\ell_1 +1 \leq q+t\ell_1$ (where the inequality holds since $\ell_2 >0$ and thus $\ell_1 \leq t-1$),
and until position $q+t(\ell_1-1)+\ell_1 + t-\ell_1  = q+t\ell_1 $.
Similarly to the previous case, $x_1$ intersects with $y_{q+i}$, for $(\ell_1 -1)t+1 \leq i \leq t\ell_1 - \ell_1 + 1$.

Finally, the last $\ell_2  \leq t-1$ columns of $Y$, all have a one in the last position,  and so intersect with $x_1$ that also has a one in the last position.
\end{itemize}

Hence, $X\cdot Y = C_{p,q}$ as claimed, and since the number of ones in each row/column of $X,Y$ is exactly $t$,
then $C_{p,q}$ is a submatrix of $A_{k,t}$ for any $k \geq p+q $.
See  Figure~\ref{fig:product} for an illustration of the construction described in the proof of Theorem~\ref{circularbigp}.
\EPF

\section{The range of large $p$}
\label{seclargep}

Bollob{\'a}s~\cite{bollobas1965generalized} proved that for any $m$ pairs of subsets
$(A_i,B_i)$, such that $|A_i| = a, |B_i| =b$ for $1 \leq i \leq m$,
and  $A_i \cap B_j = \emptyset$ if and only if $i=j$, it holds that $m  \leq {a+b \choose a}.$
An immediate corollary of this theorem is that the largest circulant submatrix $C_{p,q}$ of $A_{k,t}$,  for $q = 1$,
is of size ${2t \choose t}\times {2t \choose t}$, and this result is tight.

This theorem has several generalizations, among them is a result of
Frankl~\cite{frankl1982extremal} and Kalai~\cite{kalai1984intersection} that considered the skew version of the problem,
and showed that the same bound holds even under the following relaxed assumptions:
Let $(A_i,B_i)$ be pairs of sets, such that $|A_i| = a, |B_i| =b$ for $1 \leq i \leq m$,
$A_i \cap B_i = \emptyset$ for every $1 \leq i \leq m$, and $A_i \cap B_j \neq \emptyset$ if  $i>j$.
Then $m  \leq {a+b \choose a}$. From this we get immediately the following corollary:

\begin{coro}
\label{uppercorolarge}
Let $C_{p,q}$ be a submatrix of $A_{k,t}$ of size $n\times n$, for a given fixed $q$.
Then,   $n \leq {2t \choose t} + q-1$, that is, $p \leq {2t \choose t} -1$.
\end{coro}
\BPF
Assume that the first row of $C_{p,q}$ has $q$ zeros followed by $p$ ones
(this can be achieved simply by moving the last row of $C_{p,q}$, as defined in the introduction, to be the first row).
Now consider the submatrix $B$ of $C_{p,q}$ that is defined by the first $p+1$ rows and columns of $C_{p,q}$.
The matrix $B$ maintains the conditions of the Theorem of Frankl and Kalai, and, thus, its size is at most ${2t \choose t}\times {2t \choose t}$.
Hence,  $n \leq {2t \choose t} + q-1 $ as claimed, and $p \leq {2t \choose t}-1$.
\EPF\\

The following lemma presents a simple construction of a large circulant submatrix $C_{p,q}$ of $A_{k,t}$ for a given fixed $q$.

\BL
\label{simpleconstruction}
Let $q >0$, $t \geq q$, where  $t \bmod q = 0$. Then $C_{p,q}$ is a submatrix of $A_{k,t}$, for $p = q\cdot ({2t/q \choose t/q} -1)$ and $k \geq 3t - t/q$.
\EL
\BPF
Let $n = {2t/q \choose t/q}$.
The matrix  $C_{p,q}$, where $p + q = q\cdot n$, can be partitioned into $q$ disjoint submatrices of size $n \times (p+q)$,
as follows. The  $i$'th submatrix, $1 \leq i \leq q$,
contains rows $i+j\cdot q$, $0 \leq j \leq n-1$, of $C_{p,q}$.
Each such submatrix is a blowup of a matrix that represents the bipartite crown graph of size $n \times n$.
That is, partition each row of these $q$ submatrices into blocks of $q$ consecutive entries,
where the blocks of the $i$'th submatrix are shifted by one position compared to the blocks of the previous submatrix  (in a circulant way).
Thus, the entries in each block  are identical (either all ones or all zeros).

Hence, each one of these $q$ submatrices is the intersection matrix of the two families of all subsets of size $t/q$ of $[2t/q]$.
Now, if we take disjoint copies of these subsets, and label each column in $C_{p,q}$ with the subset that is the union of all subsets of size $t/q$ assigned to this column,
then we get $q \cdot n$ subsets of size $t$ assigned to the columns.
As to the rows, each row is assigned a different subset of size $t/q$, and therefore, we can add to each such subset
$t - t/q$ identical elements, so that the rows are also assigned subsets of size $t$.
Finally, by taking the union of all subsets, we get that $k \geq 2t + t-t/q = 3t - t/q$.
\EPF\\

The size of the construction given in Lemma~\ref{simpleconstruction} is quite far from the upper bound given in Corollary~\ref{uppercorolarge}.
As we show in the next subsection, there exists a slightly larger construction for $q = 2$.
Finally, in Subsection~\ref{larget=2}, we show that the upper bound of Corollary~\ref{uppercorolarge} is tight
for $t=2, p = {2t \choose t}- 1 = 5$ and $q =1, 3$, but there is no $k$ for which $C_{5,2}$ is a submatrix of $A_{k,2}$ when $q \neq 1,3$.

\subsection{$ q = 2$ and $p = 2^t + 2^{t-2} -2$}

We now prove Theorem~\ref{theo:largep} and show a construction of $C_{p,q}$ for $q = 2$ and $p$ that is exponential in $t$.
Throughout this section, we  assume that the first row of $C_{p,q}$ is $(0,0,1,1,...,1)$ and its last row is $(0,1,1,1...,1,0)$.

\begin{figure}[htb!]
\captionsetup{width=0.9\textwidth}
\centering
    \includegraphics[width=0.7\textwidth]{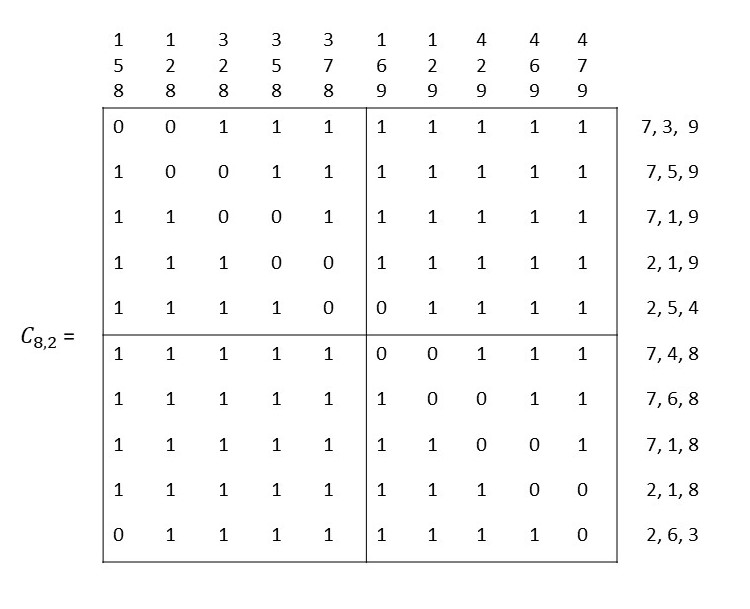}
\caption{\small The construction of $C_{p,q}$ described in Theorem~\ref{theo:largep}, for $q = 2$ and $p = 2^t + 2^{t-2} -2$, where $t = 3$.
The column indices are written above the matrix $C_{8,2}$ and the row indices to the right of the matrix.}
\label{fig:smallmatrix}
\end{figure}

The construction we present is recursive in nature, and exploits the fact that $C_{p,2}$ has two blocks on the main diagonal, such that each one of these  blocks is half the size of $C_{p,2}$,
and the structure of each block is almost identical to that of $C_{p,q}$, where the only difference is that there is a $1$ in the first position of the last row instead of a zero in $C_{p,q}$.
This small difference complicates the recursive argument. The details of the proof follow.
See Figure~\ref{fig:smallmatrix} for an illustration of the proof of Theorem~\ref{theo:largep}.

\BT
\label{theo:largep}
Let $q = 2$ and $p = 2^t + 2^{t-2} -2$, where $t > 2$. Then $C_{p,q}$ is a  submatrix of $A_{k,t}$, for large enough $k$.
\ET
\BPF
Let $h = (p+q)/2$.
We prove by induction on $t$ that $C_{p,q}$ is the intersection matrix of two  families of $t$-subsets
$$\mathcal{F}_{a,b} = \{F_1,...,F_{p+q}\}, \ \ \ \ \ \ \ \ \mathcal{G}_{a,b} = \{G_1,...,G_{p+q}\},$$
where the subsets in $\mathcal{F}_{a,b}$ are the row indices and the subsets in $\mathcal{G}_{a,b}$ are the column indices,
and $a,b$ are two integers with the following properties:
\begin{itemize}
\item
 $a \in G_1,...,G_{h}$, and  $a \in F_{h+1},...,F_{p+q-1}$.
 \item
$b \in G_{h+1},...,G_{p+q}$, and $b \in F_1,...,F_{h-1}$.
\item
 $a,b$ appear only in the subsets specified above. In particular, $a,b \not\in F_{p+q}$.
\end{itemize}
Let $\tilde{C}_{p,q}$ be the matrix that is achieved from $C_{p,q}$, by modifying to $1$ the first position of the last row of  $C_{p,q}$,
and let $\mathcal{\tilde{F}}_{a,b}$ be a family that is identical to $\mathcal{F}_{a,b}$ with one difference: the subset $F_{p+q}$ also contains  the element $a$.
It is not hard to verify that if $C_{p,q}$ is the intersection matrix of $\mathcal{F}_{a,b}$ and $\mathcal{G}_{a,b}$, then
$\tilde{C}_{p,q}$ is the intersection matrix of $\mathcal{\tilde{F}}_{a,b}$ and $\mathcal{G}_{a,b}$.

The base of the induction is $t = 3$, and the construction of $C_{8,2}$ is given in Figure~\ref{fig:smallmatrix}, where in this case $a = 8, b = 9$.
Note that if we modify the last row index $\{2,6,3\}$ to be $\{2,6,3,8\}$, then we get a construction of $\tilde{C}_{8,2}$ as claimed.

Assume now that $t > 3$, let $p_t = 2^t + 2^{t-2} -2$ and $p_{t-1} = 2^{t-1} + 2^{t-3} -2$, and consider $C_{p_t,q}$.
 Then it has the following structure: there are two matrices of the form $\tilde{C}_{p_{t-1},q}$ on the main diagonal,
and two blocks of size $(p_t+q)/2$ that are all one, but the leftmost entry on the bottom row of each of these blocks that is a $0$.

By the induction hypothesis there exist, as specified above, two families of $(t-1)$-subsets
$$\mathcal{F}_{a,b} = \{F_1,...,F_{p_{t-1}+q}\}, \ \ \ \ \ \ \ \ \mathcal{G}_{a,b} = \{G_1,...,G_{p_{t-1}+q}\},$$
whose intersection matrix is $ C_{p_{t-1},q}$.


Let $\mathcal{F'}_{b,a} = \{F'_1,...,F'_{p_{t-1}+q}\}$  be a family of subsets that is identical to $\mathcal{F}_{a,b}$,
but  $a,b$ are interchanged in all subsets. That is, for $1 \leq i \leq p_{t-1} + q$:
$$
F'_i = \left\{
  \begin{array}{ll}
     F_i \setminus \{a\} \cup \{b\}, & \hbox{if $a \in F_i$,} \\
     F_i \setminus \{b\} \cup \{a\}, & \hbox{if $b \in F_i$,} \\
     F_i, & \hbox{if $a,b \not \in F_i$.}
  \end{array}
\right.
$$
Similarly define $\mathcal{G'}_{b,a} = \{G'_1,...,G'_{p_{t-1}+q}\}$, which is identical  to $\mathcal{G}_{a,b}$,
but  $a,b$ are interchanged in all subsets.
Note that since $a,b$ appear only in subsets as specified above, then it also holds that $C_{p_{t-1},q}$ is the intersection matrix
of the two families $\mathcal{F'}_{b,a}$  and $\mathcal{G'}_{b,a}$.


Now let $c,d$ be two new elements that do not appear in any of the above families,
and define the following families:

\begin{eqnarray*}
 \mathcal{F}_{d}  &=& \{F_1 \cup \{d\},F_2 \cup \{d\},...,F_{p_{t-1}+q-1} \cup \{d\},F_{p_{t-1}+q} \cup \{a\}\},  \\
& & \\
 \mathcal{F}_{c}   & =& \{F'_1 \cup \{c\},F'_2 \cup \{c\},...,F'_{p_{t-1}+q-1} \cup \{c\},F'_{p_{t-1}+q} \cup \{b\}\},  \\
& & \\
 \mathcal{G}_{c}  &=& \{G_1 \cup \{c\},G_2 \cup \{c\},...,G_{p_{t-1}+q} \cup \{c\}\}, \\
& & \\
  \mathcal{G}_{d}   & =& \{G'_1 \cup \{d\},G'_2 \cup \{d\},...,G'_{p_{t-1}+q} \cup \{d\}\}.
\end{eqnarray*}

Finally, define the families $\mathcal{F}_{c,d}, \mathcal{G}_{c,d}$ as follows:
$$ \mathcal{F}_{c,d} = \mathcal{F}_{d} \cup \mathcal{F}_{c}, \ \ \ \ \ \ \  \mathcal{G}_{c,d}  =  \mathcal{G}_{c}   \cup  \mathcal{G}_{d}.
$$
It is clear that $\mathcal{F}_{c,d}, \mathcal{G}_{c,d}$ are two families of $t$-sets, each of size $p_t +q$,
and their structure is as claimed above, where $c$ and $d$ are in the role of $a$ and $b$, respectively.
It remains to prove that $C_{p_t,q}$ is the intersection matrix of $\mathcal{F}_{c,d}, \mathcal{G}_{c,d}$.
First note that by the induction hypothesis, and using the structure of the subsets we defined,
$\tilde{C}_{p_{t-1,q}}$ is the intersection matrix of $\mathcal{F}_{d}, \mathcal{G}_{c}$,
as well as the intersection of $\mathcal{F}_{c}, \mathcal{G}_{d}$.

Consider now the matrix $C$ which is the intersection matrix of $\mathcal{F}_{d}, \mathcal{G}_{d}$.
It is clear that the first $p_{t-1}+q-1$ rows of $C$ are all ones, since
the first $p_{t-1}+q-1$ families of  $\mathcal{F}_{d}, \mathcal{G}_{d}$ all contain $d$.
We next show that the last row of $C$ is of the form $(0,1,1,...,1)$.
By the induction hypothesis, the intersection of $F_{p_{t-1}+q}$ with all subsets of $\mathcal{G}_{a,b}$ gives a vector of the form $(0,1,1,...1,0)$.
Thus, since $a,b \not\in F_{p_{t-1}+q}$ and $\mathcal{G'}_{b,a}$ is identical  to $\mathcal{G}_{a,b}$,
but  $a,b$ are interchanged in all subsets, then the intersection of
of $F_{p_{t-1}+q}$ with all subsets of $\mathcal{G'}_{b,a}$ results also with the vector $(0,1,1,...1,0)$.
Since the last subset of $\mathcal{F}_{d}$ is defined as  $F_{p_{t-1}+q} \cup \{a\}$ and the last subset of  $\mathcal{G}_{d}$ is $G'_{p_{t-1}+q} \cup \{d\}$,
and  $a \in G'_{p_{t-1}+q}$, then we get that $F_{p_{t-1}+q} \cup \{a\}$ and $G'_{p_{t-1}+q} \cup \{d\}$ also intersect as required.

A similar argument shows that the intersection matrix of $\mathcal{F}_{d}, \mathcal{G}_{d}$
is also a matrix that is all one, but the first element on the last row of this matrix, which is a zero.
This completes the proof of the theorem. \EPF\\

\subsection{The values $t = 2, p = {2t \choose t} -1, q > 0$}
\label{larget=2}
Finally, we  address the range of values of $t = 2$ and $p = {2t \choose t} -1 = 5$.
We first show that $C_{p,q}$ is a submatrix of $A_{k,t}$ for these values of $p$ and $t$, and for $q = 1,3$.

\begin{figure}[htb!]
\captionsetup{width=0.9\textwidth}
\centering
    \includegraphics[width=0.8\textwidth]{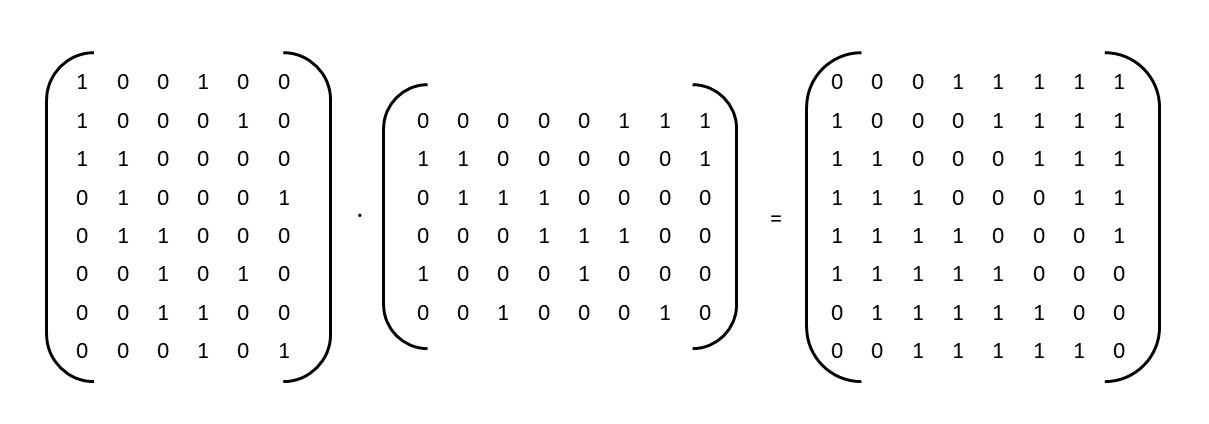}
\caption{\small A construction of $C_{p,q}$ for $t = 2, p =5, q = 3$.}
\label{fig:product1}
\end{figure}

\BL
Let $t = 2$ and $p = {2t \choose t}-1 = 5$.
Then $C_{p,q}$ is a submatrix of $A_{k,t}$ for $q = 1$ and $k \geq 5$, or for $q = 3$ and $k \geq 6$.
\EL
\BPF
If  $t = 2, p = 5, q = 1$, then $C_{5,1}$ is a submatrix of $A_{5,2}$. Simply take as row/column indices all subsets of size $2$ of $[4]$.
As to the case of  $t = 2, p = 5,  q = 3$, Figure~\ref{fig:product1} shows that $C_{5,3}$ is a submatrix of $A_{k,2}$, for any $k \geq 6$.
\EPF\\

We conclude by proving that $C_{p,q}$ is not  a submatrix of $A_{k,t}$ for $t =2, p = 5$ and $q \neq 1,3$.
Unfortunately, this proof cannot be generalized to the case of $p = {2t \choose t}-1$ and $t > 2$.
Thus, it remains an open problem to determine for what values of $q> 1$ is $C_{p,q}$ a submatrix of $A_{k,t}$, when
$p = {2t \choose t}-1$ and $t > 2$.

\BL
Let $t = 2$, $p = {2t \choose t}-1 = 5, q \neq 1,3$, $q > 0$. Then $C_{p,q}$ is not a submatrix of $A_{k,t}$ for any $k$.
\EL
\BPF
Assume by contradiction that $C_{5,q}$ is  a submatrix of $A_{k,2}$ for some $k$, where the first row of $C_{5,q}$
starts with $q \neq 1,3$ zeros, followed by $p$ ones.
Let $n = p+q = q+5 \geq 7$ be the size of  $C_{5,q}$,
and let $A_i, B_i$, $0 \leq i \leq n-1$, be the $2$-uniform subsets defining the row and column indices, respectively, of $C_{5,q}$.

Assume first that there exists some $i$ such that $B_i \cap B_{(i+1)\bmod n} = \emptyset$, that is,
two consecutive column indices are disjoint.
Since $C_{5,q}$ is circulant, then we can assume that $i = 0$, that is, $B_0 \cap B_1 = \emptyset$.
Since $B_0$ and $B_1$ both intersect with $A_3,A_4,A_5$,  then each of these three subsets
contains one element from each of $B_0,B_1$. Furthermore, as all subsets are different and of size $2$,
then each element of $B_0,B_1$ is contained in at most two of these three subsets.

Next consider $B_2$. It also intersects with $A_3,A_4,A_5$, and since there is no common element of $B_0,B_1$ in these three subsets,
then $B_2$ also includes two elements from $B_0 \cup B_1$ (although here $B_2$ can contain two elements from the same subset $B_0$ or $B_1$).

Now, consider $A_{7 \bmod n}$, where if $q = 2$ then $A_{7 \bmod n} = A_0$ and otherwise, $A_{7 \bmod n} = A_7$.
In both cases, since $p = 5$, the row labeled by $A_{7 \bmod n}$ starts with two zeros followed by a one.
Thus, since $A_{7 \bmod n} \cap B_2 \neq \emptyset$, then $A_{7 \bmod n}$ contains an element from $B_0 \cup B_1$,
in contradiction to the fact that $A_{7 \bmod n} \cap B_0 = A_{7 \bmod n} \cap B_1 = \emptyset$.

Hence, we can assume from now on that  $B_i \cap B_{(i+1)\bmod n} \neq \emptyset$, and similarly
that $A_i \cap A_{(i+1)\bmod n} \neq \emptyset$, for $0 \leq i \leq n-1$.
There are two cases:
\begin{itemize}
\item
There exists an $i$ such that $B_i \cap  B_{(i+1)\bmod n} \cap  B_{(i+2)\bmod n} \neq \emptyset$.
Since $C_{p,q}$ is circulant, then assume that $i = 0$, and let $b \in B_0,B_1,B_2$.
Thus,  $B_0=\{b_0,b\}, B_1=\{b_1,b\},B_2=\{b_2,b\}$.
From this and the structure of $C_{5,q}$, we can deduce the following:
\begin{enumerate}
\item
$b_0 \in A_1$, $A_2 = \{b_0,b_1\}$, and $b \in A_3,A_4,A_5$.
\item
The row labeled by $A_6$ starts with a zero followed by $5$ ones, and so $b \not\in A_6$. But $A_6 \cap B_1 \neq \emptyset, A_6 \cap B_2 \neq \emptyset$.
Thus, $A_6 = \{b_1,b_2\}$.
\item
$b\not\in B_3$ as $B_3 \cap A_3 = \emptyset$. But $B_3 \cap B_2 \neq \emptyset$, and therefore, $b_2 \in B_3$.
\item
$b, b_0, b_1 \not\in B_3$ as also $B_3 \cap A_2 = \emptyset$. But $B_3 \cap A_5 \neq \emptyset$ and $A_5 \cap A_6 \neq \emptyset$. Therefore, $b_2 \in A_5$.
\item
Since $b, b_0, b_1 \not\in B_3$ and $B_3 \cap A_4 \neq \emptyset$, then there exists a new element $b_3 \in B_3 \cap A_4$.
\item
$A_4 \cap B_4 = \emptyset$ and so $b \not \in B_4$. Hence, $b_2 \in B_4$ since $B_4 \cap A_5 \neq \emptyset$.
In a similar way, $b_1 \in B_5$.
\end{enumerate}
Hence, the subsets defining the first seven rows and columns of $C_{5,q}$ have the following structure so far,
where they are written to the left and above the submatrix:

$$
  \begin{array}{|c||c|c|c|c|c|c|c|}
\hline
          & b_0 & b_1 & b_2 & b_2  & b_2  &b_1 &  \\
          & b &  b &  b & b_3  &  & &   \\\hline\hline
       & 0   & 0   &    &       &   &   &   \\\hline
     b_0  & 1   & 0   & 0   &       &   &   &   \\\hline
  b_0,  b_1   & 1   & 1   & 0   &  0     &   &   &  \\\hline
 b        & 1   & 1   & 1   &  0     & 0  &   & \\\hline
   b, b_3      & 1   & 1   & 1   &  1     & 0  & 0  & \\\hline
   b, b_2     & 1   & 1   & 1   &  1     & 1  & 0  &0  \\\hline
  b_1,b_2        & 0   & 1   & 1   &  1     & 1  & 1  &0  \\\hline
  \end{array}
$$
Now, if $q \geq 4$, we already get a contradiction, since in $C_{5,q}$ it holds that $A_2 \cap B_5 = \emptyset$,
whereas here  $b_1 \in  A_2 \cap B_5$.

Therefore, assume that $q = 2$, and so all remaining entries in the submatrix above are ones.
From the structure of the submatrix and the information we have so far, we can deduce that $A_1 = \{b_0, b_3\}$ and hence $B_5 = \{b_1, b_0\}$
(since $B_5 \cap A_1 \neq \emptyset$ and $B_5 \cap A_4 = \emptyset$ and so  $b_3 \not\in B_5$).
But then since $b_0,b_1 \not\in A_0$, we get a contradiction since $A_0 \cap B_5 \neq \emptyset$.

\item
 $B_i \cap  B_{(i+1)\bmod n} \cap  B_{(i+2)\bmod n} = \emptyset$,  but $B_i \cap  B_{(i+1)\bmod n} \neq \emptyset$, for $0 \leq i \leq n-1$.
Thus, $B_i  = \{b_i, b_{(i+1)\bmod n}\}$, where some of the $b_i$'s may be identical.

If all $b_i$'s in the subsets $B_{q},B_{q+1},B_{q+2}, B_{q+3}, B_{q+4}$ are different, then $A_0$ cannot intersect with these subsets, since $|A_0| = 2$.
Hence, there exist $0 \leq i \neq j \leq 4$ such that $b_{q+i} = b_{q+j}$.
Assume, without loss of generality, that $i = 0$ (as the matrix is circulant).
Since  the intersection of every three consecutive subsets is empty, and each subset
contains two different elements, then $j \neq 1, 2$.
If $j = 3$ then $b_{q} = b_{q+3}$, and since $A_2$ does not intersect with $B_q, B_{q+1}$ then $b_{q}, b_{q+1}, b_{q+2} \not\in A_2$.
But $A_2 $ intersects with $B_{q+2} = \{b_{q+2}, b_{q+3}= b_q\}$, and we get a contradiction.
A similar contradiction is achieved if $j = 4$ when considering $A_5$.
\end{itemize}
Thus, in all cases we get a contradiction and the lemma follows.
\EPF\\

\subsection*{Acknowledgments}
This research did not receive any specific grant from funding agencies in the public, commercial, or not-for-profit sectors.

\bibliographystyle{plain}



\end{document}